\documentclass[11pt]{amsart}

\newcommand{\nw}{\bigwedge\nolimits}

\newcommand{\U}{\mathcal U}
\newcommand{\Ext}{Ext}

\newcommand{\PBW}{{PBW}}
\renewcommand{\phi}{\varphi}

\newcommand{\simto}{\overset{\sim}{\to}}

\newcommand{\DR}{{\mathrm{DR}}}

\newcommand{\op}{{\overline\partial}}
\newcommand{\ndot}{\bullet}

\newcommand{\g}{{\mathfrak g}}

\newcommand{\wQ}{\widetilde Q}

\def\O{\mathcal O}

\def\matho#1{\mathop{\mathrm{#1}}}

\DeclareMathSizes{11.1}{10}{8}{6}

\newcommand\nfrac[2]
{\dfrac{\raisebox{-1pt}{\fontsize{11.1}{10pt}\selectfont$#1$}}
       {\raisebox{ 1pt}{\fontsize{11.1}{10pt}\selectfont$#2$}}}

\newcommand{\Vect}{{\matho{Vect}}}

\newcommand{\Der}{\matho{Der}}
\newcommand{\Lie}{{\matho{Lie}}}

\newcommand{\Hom}{\matho{Hom}\nolimits}
\newcommand{\Fun}{\matho{Fun}}
\newcommand{\End}{\matho{End}}

\newcommand{\Tr}{{\matho{Tr}}}

\newcommand{\ad}{\matho{ad}}

\newcommand{\oTh}{\overline T_\hol}
\newcommand{\hol}{{\mathrm{hol}}}
\newcommand{\Hol}{{\mathrm{Hol}}}
\newcommand{\strange}{{\mathrm{strange}}}
\newcommand{\poly}{{\mathrm{poly}}}
\newcommand{\Rmn}{\R^{m|n}}
\newcommand{\pp}[2]{\nfrac{\partial#1}{\partial#2}}
\newcommand{\HKR}{{\mathrm{HKR}}}
\newcommand{\Coh}{{\mathrm{Coh}}}
\newcommand{\diag}{{\mathrm{diag}}}
\newcommand{\sheaf}{{\mathrm{sheaf}}}

\newcommand{\C}{\mathbb C}
\newcommand{\R}{\mathbb R}
\newcommand{\Z}{\mathbb Z}
\newcommand{\D}{\mathcal D}

\newtheorem*{theorem}{Theorem}

\newtheorem*{lemma}{Lemma}
\newtheorem*{corollary}{Corollary}
\newtheorem*{conjecture}{Conjecture}

\theoremstyle{remark}

\newtheorem*{example}{Example}
\newtheorem*{bexample}{Basic example}
\newtheorem*{examples}{Examples}

\theoremstyle{definition}
\newtheorem*{defin}{Definition}

\author{Boris Shoikhet}
\title
{On the Duflo formula for $L_\infty$-algebras and $Q$-manifolds}
\date{29.11.1998}
\address{IUM, 11 Bol'shoj Vlas'evskij per.,
Moscow 121002, Russia}
\email{borya@mccme.ru}

\begin{document}
\maketitle

\sloppy

\begin{abstract}

We prove a direct analogue of the classical Duflo formula
in the case of $L_\infty$-algebras. We conjecture an analogous
formula in the case of an arbitrary $Q$-manifold. When $G$
is a compact connected Lie group, the Duflo theorem for the
$Q$-manifold $(\Pi TG,d_\DR)$ is exactly the Duflo theorem for the
Lie algebra $\g=\Lie\, G$. The corresponding theorem for the
$Q$-manifold $(\Pi TM,d_\DR)$, where $M$ is an arbitrary smooth
manifold, is a generalization of the Duflo theorem for the case of smooth
manifolds. On the other hand, the Duflo theorem for the
$Q$-manifold $(\Pi\oTh M,\op)$, where $M$ is a complex manifold, is a
generalization of the M.\,Kontsevich's ``theorem on complex
manifold'' [K1], Sect.~8.4.
\end{abstract}

\section*{Contents}

\contentsline {section}{\tocsection {}{1}{The classical Duflo formula}}{}
\contentsline {section}{\tocsection {}{2}{Strong homotopy Lie algebras and $Q$-manifolds}}{}
\contentsline {section}{\tocsection {}{3}{Relationship with the Formality conjecture}}{}
\contentsline {section}{\tocsection {}{4}{Duflo formula for $L_\infty $-algebras}}{}
\contentsline {section}{\tocsection {}{5}{Duflo formula for $Q$-manifolds}}{}
\contentsline {section}{\tocsection {}{}{References}}{}

\section[The classical Duflo formula]{The classical Duflo formula [D], [K1]}

Let $\g$ be a finite-dimensional Lie algebra, $S^\ndot(\g)$ be
the symmetric algebra of the vector space~$\g$, and $U(\g)$ be the
universal enveloping algebra of the Lie algebra~$\g$. Both spaces
$S^\ndot(\g)$ and $U(\g)$ are $\g$-modules with respect to the
adjoint action; it follows from the Poincar\'e--Birkhoff--Witt
theorem that these modules are isomorphic. Therefore, the vector
spaces of invarinats $[S^\ndot(\g)]^\g$ and $[U(g)]^\g$ are
isomorphic. The Duflo theorem states that $[S^\ndot(\g)]^\g$
and $[U(\g)]^\g$ are canonically isomorphic as \emph{algebras}.
Moreover ([K1], Sect.~8.3) the algebras
$H^\ndot_\Lie(\g;S^\ndot(\g))$ and $H^\ndot_\Lie(\g;U(\g))$ are
canonically isomorphic.

Let us recall the construction of this isomorphism. First of all,
the map $\phi_\PBW\colon S^\ndot(\g)\to U(\g)$, defined as
follows
\begin{equation}
\phi_\PBW(g_1\cdot g_2\cdot\ldots\cdot
g_k)=\nfrac1{k!}\sum_{\sigma\in\Sigma_k}
g_{\sigma(1)}*g_{\sigma(2)}*\ldots*g_{\sigma(k)}
\end{equation}
is an isomorphism of the $\g$-modules. The corresponding map
$\phi_\PBW\colon[S^\ndot(\g)]^\g\simto[U(\g)]^\g$ is \emph{not}
an isomorphism of the algebras.

There exists an isomorphism of $\g$-modules $\phi_\strange\colon
S^\ndot(\g)\simto S^\ndot(\g)$ such that the composition
$$
[S^\ndot(\g)]^\g\overset{\phi_\strange}{\longrightarrow}
[S^\ndot(\g)]^\g\overset{\phi_\PBW}{\longrightarrow}[U(\g)]^\g
$$
is an isomorphism of the algebras.

The map $\phi_\strange\colon S^\ndot(\g)\simto S^\ndot(\g)$ is
defined as follows. Let us consider elements of~$\g^*$ as
derivations of the symmetric algebra~$S^\ndot(\g)$, then
elements of~$S^k(\g^*)$ are differential operators with constant
coefficients acting on~$S^\ndot(\g)$, If the Lie algebra~$\g$ is
finite-dimensional, there exists a canonical element $c_k\in
S^k(\g^*)$ for every $k\ge1$, namely, it is the trace of the
$k$-th power of the adjoint action,
$c_k=\{\,g\mapsto\Tr\ad^kg\,\}$. The elements~$c_k$ are
invariant elements in $S^k(\g^*)$. We set:
\begin{equation}
\phi_\strange=\exp\left(\sum_{k\ge1}\alpha_{2k}\cdot c_{2k}\right)
\end{equation}
where the rational numbers $\alpha_{2k}$ are defined by the
formula:
\begin{equation}
\sum_{k\ge1}\alpha_{2k}\cdot x^{2k}=\log\sqrt{\nfrac{e^{\frac
x2}-e^{-\frac x2}}x}.
\end{equation}
It is clear, that $\phi_\strange\colon S^\ndot(\g)\to
S^\ndot(\g)$ is an isomorphism of the $\g$-modules.

\begin{theorem}[{[D]}]
For any finite-dimensional Lie algebra~$\g$ the composition
$$
\phi_\PBW\circ\phi_\strange\colon[S^\ndot(\g)]^\g\to[U(\g)]^\g
$$
is an isomorphism of the algebras.

\qed
\end{theorem}

\begin{theorem}[{[K1]}, Sect.~8.3]
For any finite-dimensional Lie super-algebra~$\g$ the composition
$$
\phi_\PBW\circ\phi_\strange\colon H^\ndot_\Lie(\g;S^\ndot(\g))\to
H^\ndot(\g; U(\g)).
$$
is an isomorphism of the algebras.

\qed
\end{theorem}

In the present paper we prove the last Theorem also for
differential graded Lie algebras. Moreover, after minor
modifications the analogous statment is true also for strong
homotopy Lie algebras ($L_\infty$-algebras).
In fact, this result is a direct consequence of~[KSh].

\section{Strong homotopy Lie algebras and $Q$-manifolds}

\subsection{}

A strong homotopy Lie algebra~$\g$ ($L_\infty$-algebra) is, by
defenition, a $\Z$-graded vector space~$\g$ and an odd vector
field~$Q$ of degree~$+1$ on the space~$\g[1]$ such that
$[Q,Q]=0$.

On the other hand, it is an odd derivation~$Q$ of degree~$+1$ on
the algebra $\nw^\ndot(\g^*)$, such that $Q^2=0$ (here
$\nw^\ndot$ stands for the super-exterior algebra).

In the simplest case when $\g$ is a Lie algebra such a
differntial~$Q$ on~$\nw^\ndot(\g^*)$ is the cochain differential,
it contains only a quaratic part. When $\g$ is a DG Lie algebra,
the cochain differential on~$\nw^\ndot(\g^*)$ contain a
linear and a quadratic patrs, In other words, differential graded
Lie algebras give us  examples of $L_\infty$-algebras.

In the general case, the odd vector field~$Q$ contains also parts
of 3-rd, 4-th, $\dots$ degree.

A $Q$-manifold is a smooth  $\Z$-graded manifold~$X$ and an odd vector
field~$Q$ on~$X$ of degree~$+1$ such that $[Q,Q]=0$.
In other words, it is an odd
derivation~$Q$ of the algebra of smooth functions $C^\infty(X)$
such that $Q^2=0$. It is clear, that the case when $X$ is an $\Z$-graded
vector
space is exactly the case of $L_\infty$-algebras.

\begin{examples}

(1) Let $Y$ be a smooth (purely even) manifold, $TY$~be its
tangent bundle, and $X=T[1]Y$ be the
tangent bundle with fibers $T_x[1]$. The algebra of functions on
$T[1]Y$ is the algebra of differential forms on~$Y$, the de Rham
differential $d_\DR$ acts on this algebra. Therefore
$(T[1]Y,d_\DR)$ is an example of $Q$-manifold.

(2) Let $Y$ be a complex manifold, $X=\oTh[1]Y$, where $T_\hol Y$
is the holomorphic tangent bundle. Then $(\oTh[1]Y,\op)$ is an
example of $Q$-manifold. The corresponding complex is the
Dolbeault complex of the manifold~$Y$.
\end{examples}

\subsection{}

\begin{defin}[algebra of polyvector fields on a $Q$-manifold]

Let $(A,Q)$ be a differential graded commutative algebra, $\Der
A$ be the Lie algebra of derivations of the algebra~$A$. Then
$\Der A$ is a complex, the differential $D\colon\Der A\to(\Der
A)[1]$ is defined as the bracket with the derivation~$Q$:
$$
D(\xi)=[Q,\xi],\qquad\xi\in\Der A.
$$
It is clear that $D(f\cdot\xi)=Q(f)\cdot\xi+f\cdot D\xi$, $f\in A$,
$\xi\in\Der A$. The differential~$D$ acts by the Leibniz rule on
the (super) exterior algebra over the algebra~$A$,
$\nw_A^\ndot(\Der A)$, we denote the last DG algebra by
$T^\ndot_\poly(A,Q)$, the algebra of polyvector fields. When $X$
is a $Q$-manifold, the algebra $T^\ndot_\poly(X,Q)$ is, by
definition, the DG algebra $T_\poly^\ndot(C^\infty(X),Q)$.
\end{defin}

\begin{bexample}
Let $\g$ be a Lie algebra, $(\g[1],d_\Lie)$ be the corresponding
$Q$-manifold. The DG algebra of functions on this $Q$-manifold is
$C^\ndot_\Lie(\g;\C)$, the cochain complex of the Lie
algebra~$\g$. The DG algebra of polyvector fields
$T_\poly^\ndot(\g[1],d_\Lie)$ is the cochain complex
$C_\Lie^\ndot(\g,S^\ndot(\g))$ with the coefficient in the
symmetric algebra of~$\g$.
\end{bexample}

\subsection{}

\begin{conjecture}
Let $(A_1,Q_1)$ and $(A_2,Q_2)$ be two commutative  smooth DG
algebras, which are quasi-isomorphic. Then
$$
H^\ndot(T^\ndot_\poly(A_1,Q_1))\simeq H^\ndot
T^\ndot_\poly(A_2,Q_2)
$$
as algebras.
\end{conjecture}

(A DG algebra $(A,Q)$ is called smooth if it is the algebra of
functions on a smooth $Q$-manifold.)

\section{Relationship with the Formality conjecture [K1], [KSh]}

\subsection{}

Let $(\g[1],Q)$ be an $L_\infty$-algebra, $(\nw^\ndot g^*,Q)$ be
the correponding DG algebra of functions. In the case when $\g$ is
a Lie algebra, $Q=d_\Lie$, it is well-known result that
$$
HH^\ndot(C^\ndot_\Lie(\g;\C)\simeq HH^\ndot(U(\g))
$$
as algebras, where $HH^\ndot$ here stands for the Hochschild
cohomology. On the other hand,
$$
T^\ndot_\poly(\g[1],d_\Lie)=C^\ndot_\Lie(\g,S^\ndot(\g)).
$$
Therefore, the following theorem can be considered as a generalization of
the Duflo formula for the case of $L_\infty$-algebras.

\begin{theorem}
The algebras $H^\ndot(T^\ndot_\poly(\g[1],Q))$ and
$HH^\ndot(\nw^\ndot g^*,Q))$ are cannonically isomorphic.
\end{theorem}

We prove this theorem and construct an explicit isomorphism,
analogous to the isomorphism $\phi_\PBW\circ\phi_\strange$ from
Sect.~1, in the next Section. Here we explain, following [KSh],
Sect.~4, how the Duflo theorem itself follows from the Formality
conjecture of Maxim Kontsevich and his theorem on the
cup-products on tangent cohomology.

\begin{theorem}[Formality conjecture; proved in {[K1]}]
Let $T_\poly^\ndot(\Rmn)$ be the DG Lie algebra of polyvector
field on $\Rmn$ with zero-differential and the
Schouten--Nijenhuis bracket and let $\D^\ndot_\poly(\Rmn)$ be the
DG Lie algebra of polydifferential operators on~$\Rmn$ with the
Hochschild differential and the Gerstenhaber bracket. Then there
exists an $L_\infty$-quasi-isomorphism $\U\colon
T^\ndot_\poly\to\D^\ndot_\poly$.

\qed
\end{theorem}

(Here $\Rmn$ is a super-space $\R^m\oplus \R^n[1]$.)

More generally, one can consider an arbitrary finite-dimensional $\Z$-graded
vector space.

\subsection{}
An $L_\infty$-morphism of two $L_\infty$-algebras is a
$Q$-equivariant (may be nonlinear) map
$\phi\colon(\g_1[1],0)\to(\g_2[1],0)$.
In particular, if $\gamma\in(\g_1)^1$ is such that
$Q_1|_\gamma=0$ than $Q_2|_{\phi(\gamma)}=0$. In the case of DG
Lie algebra~$\g$ the equation $Q|_\gamma=0$ is exactly the
Maurer--Cartan equation:
\begin{equation}
\gamma\in\g^1,\qquad d\gamma+\nfrac12[\gamma,\gamma]=0.
\end{equation}
Also, $d+\ad\gamma$ defines a new differential on~$\g$. Moreover,
if $\phi\colon(g_1[1],0)\to(\g_2[1],0)$ is an $L_\infty$-morphism
of DG Lie algebras, it defines a map of the \emph{complexes}
$$
(\g_1,d_1+\ad\gamma)\to(\g_2,d_2+\ad\phi(\gamma))
$$
for each solution $\gamma\in\g_1^1$ of the Maurer--Cartan
equation. We denote the \emph{tangent complex} $(\g,d+\ad\gamma)$
by $T_\gamma(\g)$ and we denote by $\phi_\gamma\colon
T_\gamma(\g_1)\to T_{\phi(\gamma)}(\g_2)$ the corresponding
tangent map. The map $\phi_\gamma$ is a map of the complexes.

For any solution $\gamma\in T^1_\poly(\Rmn)$ of the
Maurer--Cartan equation there exists a product on the tangent
complex $T_\gamma(T^\ndot_\poly(\Rmn))$ which coincides with the
usual cup-product of the polyvector fields. On the other hand,
for any solution $\gamma\in\D_\poly^1(\Rmn)$ of the
Maurer--Cartan equation the usual product of Hochshild cochains:
\begin{equation}
(\phi\circ\psi)(a_1\otimes\ldots\otimes a_{k_1+k_2})=
\phi(a_1\otimes\ldots\otimes a_{k_1})\cdot
\psi(a_{k_1+1},\ldots,a_{k_1+k_2})
\end{equation}
(here $\phi\colon A^{\otimes k_1}\to A$, $\psi\colon A^{\otimes
k_2}\to A$ are Hochschild cochains on an algebra~$A$) defines a
product on the tangent complex $T_\gamma(\D_\poly^\ndot(\Rmn))$.

\begin{theorem}[Theorem on cup-products on tangent cohomology,
{[K1]}, Sect.~8]
Let $\U\colon T^\ndot_\poly(\Rmn)\to\D_\poly^\ndot(\Rmn)$ be the
Formality $L_\infty$-morphism, and let $\gamma\in
T^1_\poly(\Rmn)$ be a solution of the Maurer--Cartan equation.
Then the map
$$
[\U_\gamma]\colon H^\ndot(T_\gamma(T^\ndot_\poly))\to
H^\ndot(T_{\U(\gamma)}(\D^\ndot_\poly)),
$$
defined by the tangent map of the complexes:
$$
\U_\gamma\colon T_\gamma(T^\ndot_\poly)\to
T_{\U(\gamma)}(\D^\ndot_\poly),
$$
 is a morphism of the algebras.

\qed
\end{theorem}

\subsection{}
The differential graded Lie algebra $T^\ndot_\poly(\R^m)$ is
graded as follows:
$$
T^i_\poly(\R^m)=\{\,\text{$(i+1)$-polyvector fields}\,\}.
$$
In particular,
$$
T^0_\poly(\R^m)=\{\,\text{vector fields on $\R^m$}\,\}.
$$
Hovewer, every \emph{odd} vector field~$Q$ of degree~$+1$
on~$\Rmn$ lies in $T^1_\poly(\Rmn)$, and the Maurer--Cartan
equation is exactly the equation $[Q,Q]=0$.

Let $\g$ be purely even finite-dimensional Lie algebra, and let
$\displaystyle\gamma=\sum_{i,j,k}c_{ij}^k\xi_i\xi_j\pp{}{\xi_k}$
be the corresponding odd vector field on~$\g[1]$; the identity
$[\gamma,\gamma]=0$ is exactly the Jacobi identity.

Let us summarize some simplest facts on the tangent complex in this
case.

\begin{lemma}
\begin{enumerate}
\item
$T_\gamma(T^\ndot_\poly(\g[1]))=C_\Lie^\ndot(\g;S^\ndot(\g));$
\item $\displaystyle\U(\gamma)=\left\{f\mapsto\sum_{i,j,k}
c_{ij}^k\xi_i\xi_j\pp f{\xi_k}\right\}\in\D_\poly^1(\g[1]);$
\item the tangent complex
$T_{\U(\gamma)}(\D^\ndot_\poly(\g[1]))=CH^\ndot(C^\ndot_\Lie(\g;\C))$,
the Hochschild cohomological complex of the cochain complex of
the Lie algebra~$\g$.
\end{enumerate}
\qed
\end{lemma}

It follows from Theorem 3.2 that the Formality $L_\infty$-morphism
produces a map of the algebras
$$
[\U_\gamma]\colon H^\ndot_\Lie(\g;S^\ndot(\g))\to
HH^\ndot(C^\ndot_\Lie(\g;\C)).
$$

\subsubsection{}

\begin{lemma}
The map $[\U_\gamma]$ is an isomorphism \emph(of the vector
spaces\emph).
\end{lemma}

 \begin{proof}
 The statement of the Lemma follows from the homotopy
theory of $L_\infty$-algebras, see [K1],
 Sect.~4.5.1. If $\g_1,\g_2$ are two DG algebras, and
$$
\phi\colon(\g_1[1],0)\to(\g_2[1],0)
$$
is an $L_\infty$-\emph{quasi-isomorphism} between them, and a
solution $\gamma\in\g_1^1$ of the Maurer--Cartan equation is
\emph{sufficiently small}, i.e. lies in an open neighbourhood of~$0$ in
$\g_1[1]$, than the tangent map
$$
[\phi]\colon
H^\ndot T_\gamma(\g_1)\to H^\ndot T_{\phi(\gamma)}(\g_2)
$$
is an
\emph{isomorphism} of the vector spaces. In our case, we can consider
the vector field $\displaystyle\gamma_t=\sum_{i,j,k}t\cdot
c_{ij}^k\xi_i\xi_j\pp f{\xi_k}$, $t\in\C$, instead of the vector
field~$\gamma$. For sufficiently small $t$ the vector field
$\gamma_t$ lies in any open neighbourhood of~$0$ in
$T_\poly^\ndot(\g[1])$; on the other hand, if $[\U_{\gamma_t}]$
is an isomorphism for some $t\ne0$, than $[\U_\gamma]$ is also an
isomorphism.

\end{proof}

\section{Duflo formula for $L_\infty$-algebras}

We want to describe explicitly the tangent map
$$
\U_\gamma\colon C_\Lie^\ndot(\g;S^\ndot(\g))\to
CH^\ndot(C_\Lie^\ndot(\g;\C)).
$$

Let $\gamma\in T^1_\poly(\Rmn)$ be an arbitrary solution of the
Maurer--Cartan equation, which is a \emph{vector field}. The
following description of the tangent map $\U_\gamma$ was found
in~[KSh].

We fix a basis $x_1,\dots,x_{m+n}$ on $\Rmn$.

(i) The map $\phi_\HKR\colon
T^\ndot_\poly(\Rmn)\to\D_\poly^\ndot(\Rmn)$

This map is defined as follows:

if $\eta=\xi\wedge\dots\wedge\xi_k$, $\eta\in
T^\ndot_\poly(\Rmn)$, $\xi_1,\dots,\xi_k$ are vector fields, then
\begin{equation}
\phi_\HKR(\eta)(f_1\otimes\ldots\otimes
f_k)=\nfrac1{k!}\sum_{\sigma\in\Sigma_k}supersign(\sigma)\cdot
\xi_{\sigma(1)}(f_1)\cdot\ldots\cdot\xi_{\sigma(k)}(f_k).
\end{equation}
The $supersign(\sigma)$ is defined by the formula
$$
\xi_{\sigma(1)}\wedge\dots\wedge\xi_{\sigma(k)}=
supersign(\sigma)\cdot\xi_1\wedge\dots\wedge\xi_k.
$$
In particular, if $\xi_1,\dots,\xi_k$ are usual (even) vector
fields on~$\R^m$, then $supersign(\sigma)=sgn(\sigma)$; if all
the vector fields $\xi_1,\dots,\xi_k$ are odd, then
$supersign(\sigma)=1$ for any $\sigma\in\Sigma_k$.

(ii) the map $c_T$:
$$
\Vect(\Rmn)\to\Hom_{\Fun(\Rmn)}(\Vect(\Rmn)\to\Vect(\Rmn))
$$
(the ``Atiyah class'')

Let $\eta\in\Vect(\Rmn)$. We set
\begin{equation}
c_T(\eta)=\left\{\partial_{I(2)}\mapsto\sum_{I(1),I(3)}
\eta(dx^{I(1)})\cdot\partial_{I(1)}\partial_{I(2)}
\langle\gamma,dx^{I(3)}\rangle\partial_{I(3)}\right\}.
\end{equation}
Here $\partial_i=\pp{}{x_i}$ and $\langle
dx^i,\partial_j\rangle=\delta_{ij}$, and $I$ runs through all
possible maps
$$
I\colon\{1,2,3\}\mapsto\{1,2,\dots,m+n\}.
$$

The $k$-th power of the map $c_T$ is the map
$$
c_T^k\colon\Vect^{\otimes k}(\Rmn)\to\Hom^{\otimes
k}_{\Fun(\Rmn)}(\Vect,\Vect)
$$
There exists the trace map
$$
\Tr\colon\Hom^{\otimes k}_{\Fun}(\Vect,\Vect)\to\Fun,
$$
and the composition $\Tr\circ c_T^k$ is a map
\begin{equation}
\Tr\circ c_T^k=c_k\colon\Vect^{\otimes k}(\Rmn)\to\Fun(\Rmn).
\end{equation}
After the (super-) symmetrization we consider $c_k$ as an
operator
$$
c_k\colon T^\ndot_\poly(\Rmn)\to T^\ndot_\poly(\Rmn).
$$

\begin{example}
In the case of the odd field
$\displaystyle\gamma=\sum_{i,j,k}c_{ij}^k\xi_i\xi_j\pp{}{_k}$ on
the space $\g[1]$, where $\g$ is a finite-dimensional Lie
algebra, we have
$$
\langle\gamma,dx^{I(3)}\rangle=\sum_{i,j}c_{i,j}^{I(3)}\xi_I\xi_j.
$$
\end{example}

Let $\eta=\partial_l$ for some $l=1,\dots,\dim\g$. Then, by
formula~(7), the map
$$
c_T(\partial_l)=\left\{\partial_{I(2)}\mapsto\sum_{I(3)}
c_{l,I(2)}^{I(3)}\partial_{I(3)}\right\}=-\ad(\partial_l)
$$
coincides, up to a sign, with the ajoint action.

\begin{theorem}[{[KSh]}]
Let $\gamma\in T^1_\poly(\Rmn)$ be an odd vector field such that
$[\gamma,\gamma]=0$. Than the tangent map
$$
\U_\gamma\colon
T_\gamma(T_\poly^\ndot(\Rmn))\to
T_{\U(\gamma)}(\D_\poly^\ndot(\Rmn))
$$
coincides with the
composition
$$
\U_\gamma=\phi_\HKR\circ\phi_\strange,
$$
where
$$
\phi_\strange=\exp\left(\sum_{k\ge1}\alpha_{2k}c_{2k}\right).
$$
\emph(the rational numbers $\alpha_{2k}$ were defined in
Sect.~\emph{1)}

\qed
\end{theorem}

\begin{theorem}
Let $(\g[1],Q)$ be a finite-dimensional $L_\infty$-algebra, and
let $Q=Q^{(1)}+Q^{(2)}+Q^{(3)}+\dots$, where $Q^{(i)}$ is the
part of the vector field~$Q$ of $i$-th degree. Then the map
$$
\phi_\HKR\circ\phi_\strange\colon
H^\ndot(T^\ndot_\poly(\g[1],Q))\to HH^\ndot(\nw^\ndot g^*,Q))
$$
is an isomorphism of the algebras, where the operators
$$
c_{2k}\colon T^\ndot_\poly(\g[1],Q)\to T_\poly^\ndot(\g[1],Q)
$$
are defined by formulas \emph{(7)} and \emph{(8)} when
$\gamma=Q$.

\qed
\end{theorem}

\begin{corollary}
Let $\g$ be a finite-dimensional differential graded Lie algebra,
i.e. $Q=Q^{(1)}+Q^{(2)}$. Then $c_T,c_{2k}$ do not depend
on~$Q^{(1)}$, i.e. on the differential in the DG LIe
algebra~$\g$. Therefore, the usual Duflo formula (Section~1)
defines an isomorphism of the algebras
$$
\phi_\PBW\circ\phi_\strange\colon H^\ndot_\Lie(\g;S^\ndot(\g))\to
HH^\ndot(C^\ndot_\Lie(\g;\C)).
$$
\end{corollary}

\begin{proof}
The expression for the ``Atiyah class'' $c_T$ depends only on the
second derivatives of the vector field~$Q$ (see formula~(7)), but
$Q^{(1)}$ is the linear term.

\end{proof}

\section{Duflo formula for $Q$-manifolds}

Here we propose a conjecture giving an explicit formula for the
isomorphism of the algebras:
$$
H^\ndot(T^\ndot_\poly(C^\infty(X),Q))\to HH^\ndot(C^\infty(X),Q)
$$
for any smooth $Q$-manifold~$X$.

\subsection{}
Let us explain why the Duflo isomorphism for a smooth
$Q$-manifold $X$ should exist.

\begin{theorem}[Formality conjecure for smooth manifold; proved
in~{[K1]}, Sect.~7]

Let $X$ be a smooth super-manifold, $T^\ndot_\poly(X)$ be the DG
Lie algebra of smooth polyvector fields on~$X$, and let
$\D^\ndot_\poly(X)$ be the DG Lie algebra of smooth
polydifferential operators on~$X$. Then there exists an
$L_\infty$-quasi-isomorphism $\U\colon
T^\ndot_\poly(X)\to\D_\poly^\ndot(X)$.

\qed
\end{theorem}

More generally, one can consider any smooth $\Z$-graded manifold~$X$.

There does not exist a canonical choice of this
$L_\infty$-quasi-isomorphism. It was constructed in~[K1], Sect.~7
canonically ``up to a homotopy''.

\begin{conjecture}[{[K2]}, Sect.~6.4]
Let $X$ be a smooth super-manifold, $\U$ be an
$L_\infty$-quasi-isomorphism $\U\colon
T^\ndot_\poly(X)\to\D^\ndot_\poly(X)$ constructed in~\emph{[K1]},
Sect.~\emph7, and let $\gamma\in T^1_\poly(X)$ be such that
$[\gamma,\gamma]=0$ \emph(the Maurer--Cartan equation\emph). Then
the map
$$
[\U_\gamma]\colon H^\ndot T_\gamma(T^\ndot_\poly(X))\to H^\ndot
T_{\U(\gamma)}(\D^\ndot_\poly(X))
$$
induced by the tangent map~$\U_\gamma$, is a morphism of the
algebras.
\end{conjecture}

One can also suppose, that $\U_1=\phi_\HKR$ and that
$\U(\gamma)=\phi_\HKR(\gamma)$.

Let $X$ be a smooth super-manifold, $Q$ be an odd vector field of
degree~$+1$ on~$X$ such that $[Q,Q]=0$ (i.e., $X$ is a
$Q$-manifold). It follows from the above Conjecture that the
tangent map
$$
[\U_Q]\colon H^\ndot T^\ndot_\poly(C^\infty(X),Q)\to
HH^\ndot(C^\infty(X),Q)
$$
is a morphism of the algebras. Then the arguments analogous to
Lemma~3.3.1 shows that $[\U_Q]$ is in fact an
\emph{isomorphism} of the algebras.

It would be very interesting to find a description of the tangent
map
$$
\U_Q\colon T^\ndot_\poly(C^\infty(X),Q)\to
CH^\ndot(C^\infty(X),Q),
$$
analogous to the description given in Theorem~6 in the local
case. The problem is that the $L_\infty$-quasi-isomorphism
$\U\colon T^\ndot_\poly(X)\to\D^\ndot_\poly(X)$ is not defined
canonically and therefore the tangent map~$\U_Q$ also is not
defined canonically. On the other hand, the map
$$
[\U_Q]\colon H^\ndot T^\ndot_\poly(C^\infty(X),Q)\to
HH^\ndot(C^\infty(X),Q)
$$
is defined canonically. The Conjecture~5.2 below (Duflo formula
in the case of $Q$-manifolds) describes explicitly the
map~$[\U_Q]$.

\subsubsection{}

One can apply Conjecture~2.3 instead of Conjecture~5.1. Indeed,
the DG algebra $(C^\infty(X),Q)$ is quasi-isomorphic to a
finite-dimensional $L_\infty$-algebra $(\g[1],\wQ)$. It is easy
to see that the Hochschild cohomology of quasi-isomorphic DG
algebras coincides. On the other hand, Conjecture~2.3 states that
$$
H^\ndot(T^\ndot_\poly(C^\infty(X),Q))\simeq
H^\ndot(T^\ndot_\poly(\g[1],\wQ))
$$
as algebras. However, this approach does not lead us to the
explicit formula.

\subsection{}

In this section we formulate a Conjecture about the tangent
map~$[\U_Q]$.

\subsubsection{The Atiyah class in Lie algebra cohomology.}

Let $X$ be a smooth super-manifold, $\Vect(X)$ be the graded Lie
algebra of the smooth vector fields on~$X$. Let $TX$ be the
tangent bundle of the manifold~$X$. 
We denote by $J^1(TX)$ the bundle of 1-jets of the tangent
bundle. The space of global sections $\Gamma_X(J^1(TX))$ has
a natural structure of $\Vect(X)$-module.
There exists the canonical map of the $\Vect(X)$-modules
$$
\widetilde p\colon\Gamma_X(J^1(TX))\to\Vect(X).
$$
It is clear that the kernel of
the map~$\widetilde p$ is the $\Vect(X)$-module
$\Omega_X^1\otimes_{C^\infty(X)}\Vect(X)$. We obtain a short
exact sequence
\begin{equation}
0\to\Omega_X^1\otimes_{C^\infty(X)}\Vect(X)\to\Gamma_X(J^1(TX))\to
\Vect(X)\to0.
\end{equation}
Let us note that both maps in~(9) are maps of
$C^\infty(X)$-modules. The short exact sequence~(9) defines the
``Atiyah class''
\begin{equation}
c_T\in C^1_\Lie(\Vect(X);\Omega_X^1\otimes_{C^\infty(X)}
\Hom_{C^\infty(X)}(\Vect(X),\Vect(X)).
\end{equation}

\subsubsection{}

Let $X$ be a smooth $Q$-manifold.The value $c_T(Q)$ of the Atiyah
class~(10) on the odd vector field~$Q$ gives us an element
$$
c_T(Q)\in\Omega_X^1\otimes_{C^\infty(X)}\Hom_{C^\infty(X)}
(\Vect(X),\Vect(X)).
$$
This is an explicit analogoue of the Atiyah class in the case of
$L_\infty$-algebras given by formula~(7).

The $k$-th power of the element $c_T(Q)$ is a map
$$
c_T^k(Q)\colon\Vect^{\otimes k}(X)\to\Hom^{\otimes
k}_{C^\infty(X)}(\Vect(X),\Vect(X)).
$$
Futhermore, there exists the trace map
$$
\Tr\colon\Hom^{\otimes k}_{C^\infty(X)}(\Vect(X),\Vect(X))\to
C^\infty(X),
$$
and we obtain an element
\begin{equation}
c_k=\Tr\circ c_T^k(Q)\colon\Vect^{\otimes k}(X)\to C^\infty(X)
\end{equation}
After the (super-) symmetrization we can consider $c_k$ as
operators
$$
c_k\colon T^\ndot_\poly(X)\to T^\ndot_\poly(X).
$$

\begin{lemma}
The map $c_k$ is a map of the complexes
$$
c_k\colon T^\ndot_\poly(X,Q)\to T^\ndot_\poly(X,Q).
$$
\end{lemma}

\begin{proof}
It is suuficient to prove that $[Q,c_T(Q)]=0$. But $c_T$ is a
1-cocycle on $\Vect(X)$, and therefore  $(d_\Lie c_T)(Q,Q)=0$. By
the definition of the cochain differential $d_\Lie$, we have
$$
(d_\Lie c_T)(Q,Q)=c_T([Q,Q])-2[Q,c_T(Q)].
$$
The desired result follows now from the identity $[Q,Q]=0$.

\end{proof}

We set:
$$
\phi_\strange=\exp\left(\sum_{k\ge1}\alpha_{2k}c_{2k}\right)\colon
T^\ndot_\poly(X,Q)\to T^\ndot_\poly(X,Q),
$$
where the numbers~$\alpha_{2k}$ are defined by formula~(3).

\begin{conjecture}[Duflo formula for $Q$-manifolds]
Let $X$ be a smooth $Q$-manifold.

\begin{enumerate}
\item the map
$$
[\phi_\HKR\circ\phi_\strange]\colon
H^\ndot(T^\ndot_\poly(X,Q))\to HH^\ndot(C^\infty(X),Q)
$$
is an isomorphism of the algebras\emph;
\item the map $[\U_Q]$, induced by the tangent map~$\U_Q$,
coincides with the map $[\phi_\HKR\circ\phi_\strange]$.
\end{enumerate}
\end{conjecture}

\subsection{Examples}

\subsubsection{\emph(the de Rham complex\emph).}

Let $Y$ be a (purely even) smooth manifold, and let
$X=(T[1]Y,d_\DR)$. We consider $X$ as a $Q$-manifold, the DG
algebra of functions $C^\infty(X,Q)$ coincides with the de Rham
complex of the manifold~$Y$.

In the case when $Y=G$ be a (connected compact) Lie group, the
corresponding Duflo formula for $Q$-manifold $(T[1]G,d_\DR)$ can
be considered as the classical Duflo formula
for the Lie algebra~$\g$. We consider the Duflo formula
for the $Q$-manifold
$X=(T[1]Y,d_\DR)$ as an analogoue of the classical Duflo formula
for smooth manifolds.

\subsubsection{\emph(the Dolbeault complex\emph).}

Let $Y$ be a complex manifold, $T_\hol Y$ be its holomorphic
tangent bundle, $\oTh Y$ be its anti-holomorphic tangent bundle. We
consider $X=(\oTh[1]Y,\op)$ as a $Q$-manifold, the DG algebra of
functions $C^\infty(X,Q)$ coincides with the Dolbeault complex of
the manifold~$Y$.

There exist at least two different ways to define the notion of
the ``Hochschild cohomology of the structural sheaf~$\O_Y$.''

\noindent
5.3.2.1 (M.\,Kontsevich). One can define $HH^\ndot(\O_Y)$ as the
algebra of $\Ext$-s $\Ext^\ndot_{\Coh(Y\times
Y)}(\O_\diag,\O_\diag)$. The direct analogue of the
Hochschild--Kostant--Rosenberg theorem states that
\begin{equation}
\Ext^k_{\Coh(Y\times Y)}(\O_\diag,\O_\diag)=
\bigoplus\limits_{i+j=k}H^i_\sheaf(X,T^j)
\end{equation}
(here $T^j$ be the sheaf of holomorphic $j$-polyvector fields
on~$Y$). There exist  canonical products on both sides of~(12):
the Yoneda product on the Hochschild cohomology (see formula~(5))
and the product induced by the usual cup-product of polyvector
fields on the right-hand side.

The ``theorem on complex manifold'' of M.\,Kontsevich states that
both algebras are canonically isomorphic. Let us recall the
construction of this isomorphism.

Let $\alpha_T\in H^1_\sheaf(Y,\Omega_\hol^1\otimes_\O\End
T_\hol)$ be the Atiyah class of the holomorphic tangent bundle (in
the classical sense), it define the elements
$c_k=\Tr\circ\alpha_T^k\in H^k_\sheaf(Y,\Omega_\hol^k)$, which
can be considered as the Chern classes of the tangent bundle.
One can consider $c_k$ as an operator
$$
c_k\colon H^\ndot_\sheaf(Y,T^\ndot_\poly)\to
H^\ndot_\sheaf(Y,T^\ndot_\poly).
$$

\begin{theorem}[M.\,Kontsevich]
The map
$$
\phi_\HKR\circ\phi_\strange\colon
H^\ndot_\sheaf(Y,T^\ndot_\poly)\to\Ext^\ndot_{\Coh(Y\times
Y)}(\O_\diag,\O_\diag)
$$
is an isomorphism of the algebras.

\qed
\end{theorem}

\noindent
5.3.2.2. We define the Hochschild cohomology $HH^\ndot(\O_Y)$ of
the structural sheaf~$\O_Y$ as the Hochschild cohomology of the
corresponding Dolbeault complex,
$HH^\ndot(C^\infty(X),\op)$. We claim, that this definition does \emph{not}
coincide with the definition given in Sect.~5.3.2.1.
Conjecture~5.2 states, that the algebras
$HH^\ndot(C^\infty(X),\op)$ and $H^\ndot
T^\ndot_\poly(C^\infty(X),\op)$ are isomorphic, and gives an
explicit formula for the isomorphism.

We claim, that $\bigoplus\limits_{i,j}H^i(Y,T^j_\poly)$ is a
\emph{proper} subalgebra of the algebra
$H^\ndot T^\ndot_\poly(C^\infty(X),\op)$, and
$\Ext^\ndot_{\Coh(Y,Y)}(\O_\diag,\O_\diag)$ is a \emph{proper}
subalgebra of the algebra $HH^\ndot(C^\infty(X),\op)$. Indeed,
let us consider the derivations of the Dolbeault complex of the
manifold~$Y$ of the form $C^\infty(X)\otimes_\O w$, where $w$ is
a \emph{holomorphic} vector field on~$Y$; the differential
$\ad(\op)$ on these derivations is the same as the
differntial in the Dolbeault
complex $D^\ndot(\Hol)$ of the sheaf of holomorphic vector fields
on~$Y$. It is clear that
$\nw^\ndot_{(C^\infty(X),\op)}D^\ndot(\Hol)$ is a proper DG
subalgebra in $T^\ndot_\poly(C^\infty(X),\op)$;
on the other hand,
$$
H^\ndot\left(\nw^\ndot_{(C^\infty(X),\op)}D^\ndot(\Hol)\right)\simeq\bigoplus
H^\ndot_\sheaf(Y,T^\ndot_\poly).
$$

Therefore, Conjecture 5.2 in the case of the $Q$-manifold
$X=(\oTh[1]Y,\op)$ gives us a generalization of the
M.\,Kontsevich's theorem on complex manifold.

\section*{Acknowledgements}

I am grateful to Boris Feigin and to Maxim Kontsevich for many
useful discussions.

\end{document}